\documentclass[a4paper]{article}
\pdfoutput=1
\usepackage{graphicx,wrapfig}
\usepackage{lipsum}
\usepackage{amsmath,amssymb,mathrsfs}
\usepackage{hyperref}
\usepackage{url}
\usepackage{mathtools}
\usepackage{changepage}
\usepackage{multirow}
\usepackage{wrapfig}
\usepackage{subfig}
\usepackage{xcolor}
\usepackage{epsfig,enumerate,amsmath,amsfonts,amssymb,amsthm,mathrsfs,ifpdf}
\usepackage{indentfirst,relsize}
\usepackage{setspace,graphicx}
\usepackage{latexsym}
\usepackage[all]{xy}
\usepackage[usenames,dvipsnames]{pstricks}
\usepackage{pst-grad} 
\usepackage{pst-plot} 
\usepackage[margin = 2.50cm]{geometry}

\usepackage[linesnumbered,ruled,vlined]{algorithm2e}

\usepackage{color,soul}

\usepackage[colorinlistoftodos]{todonotes}

\usepackage{authblk} 

\usepackage{orcidlink}

\usepackage{enumitem}

\newcommand{\remove}[1]{}

\newtheorem{theorem}{Theorem}
\newtheorem{lemma}[theorem]{Lemma}

\newcounter{Case}[theorem]
\newtheorem{case}[Case]{Case}

\newcounter{Ca}[theorem]

\newcounter{NCa}[theorem]

\newtheorem{definition}[theorem]{Definition}

\newtheorem{observation}[theorem]{Observation}

\newtheorem{question}[theorem]{Question}

\usepackage{xcolor}
\allowdisplaybreaks

\title{Edge open packing on subclasses of chordal graphs}
\author[]{Kamal Santra \orcidlink{0009-0006-5997-1452} \footnote{kamal.7.2013@gmail.com, kamal.santra@iitg.ac.in}}
\affil[]{Department of Mathematics\\
	
	Indian Institute of Technology Guwahati\\
	
	Guwahati, 781039, Assam, India}
\date{}

\begin{document}

\maketitle
\begin{abstract}
	Packing problems in graphs are fundamental in combinatorial optimization and arise naturally in applications such as resource allocation, scheduling, and communication networks. A classical example is the \emph{induced matching} problem, where one seeks a set of edges whose endpoints induce a matching. In 2022, Chelladurai et al. introduced the notion of \emph{edge open packing}, which can be viewed as a relaxation of induced matching: instead of forcing the selected edges to induce only isolated edges, edge open packing allows them to induce disjoint stars.
	
	For a graph \(G=(V,E)\), two edges \(e_1,e_2\in E(G)\) are said to have a common edge if there exists an edge \(e\in E(G)\setminus\{e_1,e_2\}\) joining an endpoint of \(e_1\) to an endpoint of \(e_2\). A set \(D\subseteq E(G)\) is an \emph{edge open packing set} if no two edges in \(D\) have a common edge, and the maximum cardinality of such a set is the \emph{edge open packing number} \(\rho_e^o(G)\). The corresponding optimization problem is the \textsc{Maximum Edge Open Packing Problem}.
	
	In this paper, we study the computational complexity of the \textsc{Maximum Edge Open Packing Problem}. Motivated by an open question posed by Bre\v{s}ar and Samadi concerning chordal graphs, we investigate the problem on three subclasses of chordal graphs. We give an \(O(n^2)\)-time algorithm for proper interval graphs, an \(O(n+m)\)-time algorithm for block graphs, where \(n=|V(G)|\) and \(m=|E(G)|\), and an \(O(n^3)\)-time algorithm for split graphs. These results provide partial answers to the open question and contribute to the algorithmic understanding of edge packing parameters in chordal graph classes.
\end{abstract}

{\bf Keywords.}
Edge open packing, Interval graphs, Block graphs, Split graphs, Polynomial-time algorithm.

\section{Introduction}\label{EOP_Sec_1}

Throughout this paper, all graphs are finite and simple. For a graph \(G\), we denote its vertex set and edge set by \(V(G)\) and \(E(G)\), respectively. Given a vertex subset \(S\subseteq V(G)\), the induced subgraph \(G[S]\) is the subgraph with vertex set \(S\) and all edges of \(G\) whose endpoints both lie in \(S\). Similarly, for an edge subset \(D\subseteq E(G)\), the graph \(G[D]\) denotes the subgraph of \(G\) induced by the endpoints of the edges in \(D\). We also use \(G\langle D\rangle\) for the edge-induced subgraph whose edge set is exactly \(D\) and whose vertex set consists of the endpoints of the edges in \(D\). The open and closed neighbourhoods of a vertex \(v\) are denoted by \(N_G(v)\) and \(N_G[v]\), respectively. The minimum and maximum degrees of \(G\) are denoted by \(\delta(G)\) and \(\Delta(G)\). A graph is called \emph{bipartite} if its vertex set can be partitioned into two disjoint independent sets. A \emph{star} is the complete bipartite graph \(K_{1,t}\) for some \(t\geq 1\).

Packing problems in graphs form an important class of combinatorial optimization problems. They capture the task of selecting objects that are mutually compatible under some prescribed structural restriction. Classical examples include matchings, induced matchings, and several vertex and edge packing parameters. These problems arise naturally in applications such as communication networks, scheduling, resource allocation, and frequency assignment.

One important edge packing problem is the \emph{induced matching} problem. An induced matching \(M\) of a graph \(G\) is a set of edges such that no two edges of \(M\) share an endpoint and no edge of \(G\) joins an endpoint of one edge of \(M\) to an endpoint of another edge of \(M\). Equivalently, the subgraph induced by the endpoints of \(M\) is a disjoint union of \(K_{1,1}\)-stars. The maximum cardinality of an induced matching is called the induced matching number of \(G\). The problem of finding a maximum induced matching was introduced by Stockmeyer and Vazirani as the ``risk-free marriage problem''~\cite{stockmeyer1982np}, and it has since been studied extensively; see, for example,~\cite{cameron1989induced,dabrowski2013new,lozin2002maximum}.

In 2022, Chelladurai et al.~\cite{chelladurai2022edge} introduced the notion of \emph{edge open packing}. Let \(e_1,e_2\in E(G)\) be two distinct edges. We say that \(e_1\) and \(e_2\) have a \emph{common edge} if there exists an edge \(e\in E(G)\setminus\{e_1,e_2\}\) joining an endpoint of \(e_1\) to an endpoint of \(e_2\). A set \(D\subseteq E(G)\) is called an \emph{edge open packing set}, or simply an \emph{EOP set}, if no two edges of \(D\) have a common edge. The maximum cardinality of an EOP set in \(G\) is called the \emph{edge open packing number} of \(G\), and is denoted by \(\rho_e^o(G)\). An EOP set of cardinality \(\rho_e^o(G)\) is called a \(\rho_e^o(G)\)-set.

The notion of edge open packing can be viewed as a natural relaxation of induced matching. Indeed, every induced matching is an EOP set, but an EOP set need not be a matching. The selected edges in an EOP set are allowed to share a common endpoint, and therefore the subgraph \(G[D]\) induced by the endpoints of an EOP set may be a disjoint union of induced stars rather than only isolated edges. Thus, edge open packing keeps a meaningful separation condition between selected edges, while allowing a richer structure than induced matching.

Edge open packing is also closely related to injective edge coloring. An \emph{injective \(k\)-edge coloring} of a graph \(G\) is an assignment of colors from \(\{1,\ldots,k\}\) to the edges of \(G\) such that any two edges having a common edge receive distinct colors. The minimum such integer \(k\) is the \emph{injective chromatic index} of \(G\), denoted by \(\chi_i'(G)\). This notion was introduced by Cardoso et al.~\cite{cardoso2019injective} and has been further studied in~\cite{foucaud2021complexity,miao2022note}. The color classes of an injective edge coloring are precisely EOP sets. Hence, injective edge coloring can be interpreted as the problem of partitioning the edge set of a graph into EOP sets.

Chelladurai et al.~\cite{chelladurai2022edge} studied several structural aspects of edge open packing. They obtained bounds involving the diameter, size, minimum degree, clique number, and girth of a graph, and they characterized graphs \(G\) with \(\rho_e^o(G)\in\{1,2,m-2,m-1,m\}\), where \(m=|E(G)|\). Later, Bre\v{s}ar and Samadi~\cite{brevsar2024edge} studied the problem from an algorithmic point of view. They proved that the decision version of the \textsc{Maximum Edge Open Packing Problem} is NP-complete for graphs with a universal vertex, Eulerian bipartite graphs, and planar graphs of maximum degree at most \(4\). On the positive side, they gave a linear-time algorithm for trees. They also characterized graphs attaining the bound \(\rho_e^o(G)\leq |E(G)|/\delta(G)\), and studied how \(\rho_e^o(G)\) changes after deleting an edge. Further related results were obtained in~\cite{BRESAR2025,pandey2025edge}. In particular, Bre\v{s}ar et al.~\cite{BRESAR2025} studied relations between the induced matching number and the EOP number in trees and graph products, while Pandey and Santra~\cite{pandey2025edge} characterized graphs with prescribed EOP number \(t\geq 3\) and described graph families with \(\rho_e^o(G)=m-3\).

Bre\v{s}ar and Samadi~\cite{brevsar2024edge} posed the following natural question concerning chordal graphs.

\begin{question}
	Is there an efficient algorithm to determine the edge open packing number of a chordal graph?
\end{question}

In this paper, we give partial answers to this question by considering three important subclasses of chordal graphs: proper interval graphs, block graphs, and split graphs. These classes have well-known structural decompositions, and our algorithms exploit these structures. For proper interval graphs, we use a bi-compatible elimination ordering and compute the optimum over suffix subgraphs. For block graphs, we use the cut-tree and perform a bottom-up dynamic program over the rooted cut-tree. For split graphs, the structure of the clique and the independent set allows a direct characterization of the optimum value.

We recall the relevant graph classes. An edge between two non-consecutive vertices of a cycle is called a \emph{chord}. A graph is \emph{chordal} if every cycle of length at least four has a chord. A vertex \(v\) is \emph{simplicial} if \(N_G[v]\) induces a clique. A \emph{perfect elimination ordering}, or \emph{PEO}, of a graph \(G\) is an ordering \(\sigma=(v_1,v_2,\ldots,v_n)\) such that \(v_i\) is simplicial in \(G[\{v_i,v_{i+1},\ldots,v_n\}]\) for every \(i\). The existence of \emph{PEO} characterizes chordal graphs; that is, a graph $G$ is chordal if and only if $G$ has a \emph{PEO}~\cite{fulkerson1965incidence}.

A graph is an \emph{interval graph} if it is the intersection graph of a family of intervals on the real line. If no interval properly contains another, then the graph is a \emph{proper interval graph}. A PEO \(\sigma=(v_1,v_2,\ldots,v_n)\) of a chordal graph is called a \emph{bi-compatible elimination ordering}, or \emph{BCO}, if the reverse ordering is also a PEO. A graph admits a BCO if and only if it is a proper interval graph~\cite{jamison1982elimination}. We use the following standard property of BCOs.

\begin{observation}[\cite{panda2003linear}]\label{BCO_obs}
	If \(\sigma=(v_1,v_2,\ldots,v_n)\) is a \emph{\textup{BCO}} of a connected proper interval graph \(G=(V,E)\), then the set \(\{v_i,v_{i+1},\ldots,v_j\}\) induces a clique in \(G\) whenever \(v_iv_j\in E(G)\) and \(i<j\).
\end{observation}

A vertex \(v\in V(G)\) is a \emph{cut vertex} if deleting \(v\) increases the number of connected components. A \emph{block} is a maximal connected subgraph with no cut vertex. A graph is a \emph{block graph} if every block is a clique. Finally, a graph \(G=(V,E)\) is a \emph{split graph} if \(V\) can be partitioned into an independent set \(S\) and a clique \(K\).

\paragraph{\textbf{Our contributions.}}
The main results of this paper are the following.
\begin{itemize}
	\item We give an \(O(n^2)\)-time algorithm for computing the edge open packing number of a proper interval graph on \(n\) vertices.
	\item We give an \(O(n+m)\)-time algorithm for computing the edge open packing number of a block graph, where \(n=|V(G)|\) and \(m=|E(G)|\).
	\item We show that the edge open packing number of a split graph on \(n\) vertices can be computed in \(O(n^3)\) time.
\end{itemize}

The rest of the paper is organized as follows. In Section~\ref{EOP_Sec_2}, we study proper interval graphs and present an \(O(n^2)\)-time algorithm. In Section~\ref{EOP_Sec_3}, we give a dynamic-programming algorithm for block graphs running in \(O(n+m)\) time. In Section~\ref{EOP_Sec_4}, we deal with split graphs and obtain an \(O(n^3)\)-time algorithm. Finally, Section~\ref{EOP_Sec_5} concludes the paper.

\section{EOP set in proper interval graphs}\label{EOP_Sec_2}

In this section, we present a polynomial-time algorithm for computing a maximum edge open packing set in a proper interval graph. Let $G$ be a connected proper interval graph with a \emph{\textup{BCO}} $\sigma=(v_1,v_2,\ldots,v_n)$. For each vertex $v_i$, where $i\in[n]$, we define $r(v_i)$ to be the index of the rightmost vertex in the ordering $\sigma$ that is adjacent to $v_i$, including $v_i$ itself. Equivalently, $r(v_i)=\max\bigl(\{i\}\cup\{q>i: v_iv_q\in E(G)\}\bigr)$. Thus, $v_{r(v_i)}$ is the rightmost neighbour of $v_i$ in $\sigma$. Let $r(v_1)=k$. For a vertex $v_j$, let $r(v_j)=t$ be the index of the rightmost neighbour in the ordering $\sigma$, where $t\geq j$. We define $\mathcal{G}_j$ to be the subgraph of $G$ induced by all vertices appearing after $v_{r(v_j)}$ in $\sigma$; that is, $\mathcal{G}_j=G[\{v_{t+1},v_{t+2},\ldots,v_n\}]$. If $t=n$, then $\mathcal{G}_j$ is the empty graph. An edge $e = uv$ saturates its endpoints $u$ and $v$. Given a set of edges $S$, let \(V_S\) denote the set of all vertices saturated by the edges in \(S\), that is, \(V_S = V(G\langle S\rangle)\). To design the algorithm, we need the following versions of optimal EOP sets in a proper interval graph.

\begin{itemize}
	\item $\rho_{e}^{o}\left(G\right)=\max \left\{|S|: S\right.$ is an EOP set in $\left. G\right\}$;

	\item $\rho_1\left(G\right)=\max \left\{|S|: S\right.$ is an EOP set in $G$, where $v_1v_j\in S$ for some $\left. 1<j\leq k\right\}$;
	
	\begin{itemize}
		\item $\rho'\left(G\right)=\max \left\{|S|: S\right.$ is an EOP set in $G$, where $v_1v_j\in S$ for some $1<j\leq k$ and $v_jv_p \notin S$, for all $\left. k< p\leq r(v_j)\right\}$;
		
		\item $\rho''\left(G\right)=\max \left\{|S|: S\right.$ is an EOP set in $G$, where $v_1v_j\in S$ for some $1<j\leq k$ and $v_jv_p \in S$, for some $\left. k< p\leq r(v_j)\right\}$;
		
		(The above two parameters are used to calculate the value of $\rho_1\left(G\right)$)
	\end{itemize}

	\item $\rho_2\left(G\right)=\max \left\{|S|: S\right.$ is an EOP set in $G$, where $v_1v_j \notin S$, for all $\left. 1<j\leq k\right\}$;
	
\end{itemize}

Let $G$ be a connected proper interval graph with a \emph\textup{BCO} $\sigma=\{v_1, v_2, \ldots, v_n\}$ and $S$ be an arbitrary EOP set of $G$. Since the vertices of $A=\{v_1, v_2, \ldots, v_k\}$ form a clique in $G$, in $S$ at most one edge from the clique $G[A]$ can be present. Then one of the following possibilities occurs: $(i)$ $v_1v_j\in S$ for some $1<j\leq k$, or $(ii)$ $v_1$ is not incident with an edge from $S$, that is, $v_1\notin V_S$. Based on this, we have the following observation.

\begin{observation}
	Let G be a connected proper interval graph with a \emph\textup{BCO} $\sigma=\{v_1, v_2, \ldots, v_n\}$. Then $\rho_{e}^{o}\left(G\right)=\max \left\{\rho_1\left(G\right), \rho_2\left(G\right)\right\}$.
\end{observation}

In the subsequent lemmas, we determine the values of the two parameters $\rho_1\left(G\right)$, and $\rho_2\left(G\right)$, respectively.

\begin{lemma}\label{PIG_Lemma_1}
	Let $G$ be a connected proper interval graph with a \emph{\textup{BCO}} $\sigma=(v_1,v_2,\ldots,v_n)$. Also, let $r(v_1)=k$ and $\mathcal{G}_t=G[\{v_{r(v_t)+1},v_{r(v_t)+2},\ldots,v_n\}]$ for a vertex $v_t$. Then
	\begin{enumerate}[label=\textup{(\roman*)}]
		\item $\rho'(G)=1+\displaystyle\max_{1<j\leq k}\rho_e^o(\mathcal{G}_j)$;
		
		\item $\rho''(G)=2+\displaystyle\max_{\substack{1<j\leq k\\ k<p\leq r(v_j)}}\rho_e^o(\mathcal{G}_p)$.
	\end{enumerate}
	Moreover, $\rho_1(G)=\max\{\rho'(G),\rho''(G)\}$.
\end{lemma}

\begin{proof}
	Let $S$ be an EOP set of $G$ satisfying the conditions in the definition of $\rho_1(G)$. Thus $v_1v_j\in S$ for some $1<j\leq k$, where $k=r(v_1)$. Since $G$ is a proper interval graph and $\sigma$ is a BCO, by Observation~\ref{BCO_obs}, the set $\{v_1,v_2,\ldots,v_k\}$ induces a clique in $G$. For fixed $j$, since $S$ is an EOP set containing the edge $v_1v_j$, no other edge of $S$ can have an endpoint in $\{v_1,v_2,\ldots,v_{j-1},v_{j+1},\ldots,v_k\}$. In particular, no edge of the form $v_jv_r$, where $2<r\leq k$ and $r\neq j$, can belong to $S$. Moreover, at most one edge of the form $v_jv_p$, where $k<p\leq r(v_j)$, can belong to $S$, since $\{v_j,v_{j+1},\ldots,v_{r(v_j)}\}$ induces a clique in $G$. Accordingly, we consider the following two cases.

	\begin{case}\label{PIG_CA_1}
		$v_1v_j\in S$ for some $1<j\leq k$, and $v_jv_p\notin S$ for every $k<p\leq r(v_j)$.
	\end{case}
	
	Assume that $\rho'(G)=|S|$. Fix such an index $j$, and let $S_j=S\setminus\{v_1v_j\}$. By the defining condition of $\rho'(G)$, no edge of $S$ has the form $v_jv_p$ with $k<p\leq r(v_j)$. Hence, every edge of $S_j$ lies entirely in the subgraph induced by the vertices appearing after $v_{r(v_j)}$ in the ordering $\sigma$. By the definition of $\mathcal{G}_j$, this subgraph is $\mathcal{G}_j=G[\{v_{r(v_j)+1},v_{r(v_j)+2},\ldots,v_n\}]$. Therefore, $S_j$ is an EOP set of $\mathcal{G}_j$, and so $|S_j|\leq \rho_e^o(\mathcal{G}_j)$. Consequently, $\rho'(G)=|S|=1+|S_j|\leq 1+\rho_e^o(\mathcal{G}_j)\leq 1+\displaystyle\max_{1<j\leq k}\rho_e^o(\mathcal{G}_j)$.
	
	Conversely, fix an index $j$ with $1<j\leq k$, and let $T_j$ be a $\rho_e^o(\mathcal{G}_j)$-set. Define $S'=T_j\cup\{v_1v_j\}$. By the definition of $\mathcal{G}_j$, every edge of $T_j$ lies entirely in the subgraph induced by the vertices appearing after $v_{r(v_j)}$ in the ordering $\sigma$. Hence, no edge of $T_j$ has a common edge with $v_1v_j$, and so $S'$ is an EOP set of $G$. Moreover, $S'$ satisfies the defining condition of $\rho'(G)$, namely $v_1v_j\in S'$ and $v_jv_p\notin S'$ for every $k<p\leq r(v_j)$. Therefore, $\rho'(G)\geq |S'|=1+|T_j|=1+\rho_e^o(\mathcal{G}_j)$.
	
	Now choose an index $j^*$ such that $\rho_e^o(\mathcal{G}_{j^*})=\displaystyle\max_{1<j\leq k}\rho_e^o(\mathcal{G}_j)$. Then $\rho'(G)\geq 1+\rho_e^o(\mathcal{G}_{j^*})=1+\displaystyle\max_{1<j\leq k}\rho_e^o(\mathcal{G}_j)$. Combining this with the previous upper bound, we obtain $\rho'(G)=1+\displaystyle\max_{1<j\leq k}\rho_e^o(\mathcal{G}_j)$.
	
	\begin{case}\label{PIG_CA_2}
		$v_1v_j\in S$ for some $1<j\leq k$, and $v_jv_p\in S$ for some $k<p\leq r(v_j)$.
	\end{case}
	
	Assume that $\rho''(G)=|S|$. Fix such indices $j$ and $p$, and let $S_{j,p}=S\setminus\{v_1v_j,v_jv_p\}$. Since $S$ is an EOP set containing the edges $v_1v_j$ and $v_jv_p$, no edge of $S_{j,p}$ can have an endpoint among the vertices $v_1,v_2,\ldots,v_{r(v_p)}$. Indeed, such an edge would create a common edge with either $v_1v_j$ or $v_jv_p$, using the clique structure given by the BCO. Hence, every edge of $S_{j,p}$ lies entirely in the subgraph induced by the vertices appearing after $v_{r(v_p)}$ in $\sigma$. By the definition of $\mathcal{G}_p$, this subgraph is $\mathcal{G}_p=G[\{v_{r(v_p)+1},v_{r(v_p)+2},\ldots,v_n\}]$. Therefore, $S_{j,p}$ is an EOP set of $\mathcal{G}_p$, and so $|S_{j,p}|\leq \rho_e^o(\mathcal{G}_p)$. Consequently, $\rho''(G)=|S|=2+|S_{j,p}|\leq 2+\rho_e^o(\mathcal{G}_p)\leq 2+\displaystyle\max_{\substack{1<j\leq k\\ k<p\leq r(v_j)}}\rho_e^o(\mathcal{G}_p)$.
	
	Conversely, fix indices $j$ and $p$ with $1<j\leq k$ and $k<p\leq r(v_j)$, and let $T_p$ be a $\rho_e^o(\mathcal{G}_p)$-set. Define $S'=T_p\cup\{v_1v_j,v_jv_p\}$. By the definition of $\mathcal{G}_p$, every edge of $T_p$ lies after $v_{r(v_p)}$ in the ordering $\sigma$. Hence, no edge of $T_p$ has a common edge with either $v_1v_j$ or $v_jv_p$, and so $S'$ is an EOP set of $G$. Moreover, $S'$ satisfies the defining condition of $\rho''(G)$, namely $v_1v_j\in S'$ and $v_jv_p\in S'$ for some $1<j\leq k$ and $k<p\leq r(v_j)$. Therefore, $\rho''(G)\geq |S'|=2+|T_p|=2+\rho_e^o(\mathcal{G}_p)$.
	
	Now choose indices $j^*$ and $p^*$ such that $\rho_e^o(\mathcal{G}_{p^*})=\displaystyle\max_{\substack{1<j\leq k\\ k<p\leq r(v_j)}}\rho_e^o(\mathcal{G}_p)$. Then $\rho''(G)\geq 2+\rho_e^o(\mathcal{G}_{p^*})=2+\displaystyle\max_{\substack{1<j\leq k\\ k<p\leq r(v_j)}}\rho_e^o(\mathcal{G}_p)$. Combining this with the upper bound, we obtain $\rho''(G)=2+\displaystyle\max_{\substack{1<j\leq k\\ k<p\leq r(v_j)}}\rho_e^o(\mathcal{G}_p)$. This proves \textup{(ii)}.
	
	From the definitions of $\rho_1(G)$, $\rho'(G)$, and $\rho''(G)$, it follows that $\rho_1(G)=\max\{\rho'(G),\rho''(G)\}$.
\end{proof}

\begin{lemma}\label{PIG_Lemma_2}
	Let $G$ be a connected proper interval graph with a \emph{\textup{BCO}} $\sigma=(v_1,v_2,\ldots,v_n)$. For each $i\in[n]$, let $G_i=G[\{v_i,v_{i+1},\ldots,v_n\}]$. Then $\rho_2(G)=\rho_e^o(G_2)$.
\end{lemma}

\begin{proof}
	Let $S$ be an EOP set of $G$ satisfying the defining condition of $\rho_2(G)$, namely $v_1\notin V_S$, where $V_S=V(G\langle S\rangle)$ is the set of vertices saturated by the edges in $S$. Since $v_1\notin V_S$, no edge of $S$ is incident with $v_1$. Hence, every edge of $S$ lies in $G_2=G[\{v_2,v_3,\ldots,v_n\}]$, and therefore $S$ is an EOP set of $G_2$. This gives $|S|\leq \rho_e^o(G_2)$, and hence $\rho_2(G)\leq \rho_e^o(G_2)$.
	
	Conversely, let $T$ be a $\rho_e^o(G_2)$-set. Since $T$ contains no edge incident with $v_1$, we have $v_1\notin V_T$. Moreover, $T$ remains an EOP set when considered as an edge set of $G$. Indeed, any common edge between two edges of $T$ would have both its relevant endpoints in $G_2$, and hence would already be present in $G_2$, contradicting that $T$ is an EOP set of $G_2$. Thus, $T$ satisfies the defining condition of $\rho_2(G)$, and so $\rho_2(G)\geq |T|=\rho_e^o(G_2)$. Therefore, $\rho_2(G)=\rho_e^o(G_2)$.
\end{proof}

Using the structural properties of proper interval graphs and the preceding lemmas, we now give a dynamic programming algorithm for computing $\rho_e^o(G)$. The algorithm computes a bi-compatible elimination ordering only once and then works on suffix subgraphs of this ordering. For $i\in[n]$, write $r_i=r(v_i)$, and let $G_i=G[\{v_i,v_{i+1},\ldots,v_n\}]$. We store in $F[i]$ the value $\rho_e^o(G_i)$, and set $F[n+1]=0$ for the empty graph.

The dynamic program directly implements Lemmas~\ref{PIG_Lemma_1} and~\ref{PIG_Lemma_2}. For the suffix graph $G_i$, the first vertex is $v_i$. By Lemma~\ref{PIG_Lemma_2}, if $v_i$ is not saturated, then the remaining problem is exactly $G_{i+1}$, giving the candidate value $F[i+1]$. By Lemma~\ref{PIG_Lemma_1}, if $v_i$ is saturated, then either one edge $v_iv_j$ is chosen, where $i<j\leq r_i$, giving the candidate value $1+F[r_j+1]$, or two edges $v_iv_j$ and $v_jv_p$ are chosen, where $i<j\leq r_i$ and $r_i<p\leq r_j$, giving the candidate value $2+F[r_p+1]$.

Thus, $F[i]$ is obtained by taking the maximum over these possibilities. Using the monotonicity of rightmost neighbours in a BCO, the two-edge case can be simplified by taking $j=r_i$ and considering all $p$ with $r_i<p\leq r_{r_i}$. Hence, the recurrence used by the algorithm is $F[i]=\max\{F[i+1],\,1+\max_{i<j\leq r_i}F[r_j+1],\,2+\max_{r_i<p\leq r_{r_i}}F[r_p+1]\}$, where a maximum over an empty set is ignored. Accordingly, the algorithm first initializes $\mathrm{best}$ with $F[i+1]$, then updates it using the one-edge choices from Lemma~\ref{PIG_Lemma_1}(i), and finally updates it using the two-edge choices from Lemma~\ref{PIG_Lemma_1}(ii).

\begin{algorithm}[htbp!]
	\caption{\textsc{MaxEOP\_PIG}$(G)$}
	\label{EOP_PIG_Algo}
	\KwIn{A connected proper interval graph $G=(V,E)$ with $n$ vertices.}
	\KwOut{$\rho_e^o(G)$, the edge open packing number of $G$.}
	
	Compute a \emph{\textup{BCO}} $\sigma=(v_1,v_2,\ldots,v_n)$ of $G$\;
	
	\For{$i\leftarrow 1$ \KwTo $n$}{
		$r_i\leftarrow \max\bigl(\{i\}\cup\{q>i: v_iv_q\in E(G)\}\bigr)$\tcp*{$v_{r_i}$ is the rightmost neighbour of $v_i$}
	}
	
	$F[n+1]\leftarrow 0$\tcp*{Value of the empty suffix}
	
	\For{$i\leftarrow n,n-1,\ldots,1$}{
		$k\leftarrow r_i$\tcp*{$v_k$ is the rightmost neighbour of $v_i$}
		
		$\mathrm{best}\leftarrow F[i+1]$\tcp*{Case where $v_i$ is not saturated}
		
		\If{$i<k$}{
			\For{$j\leftarrow i+1$ \KwTo $k$}{
				$\mathrm{best}\leftarrow \max\{\mathrm{best},\,1+F[r_j+1]\}$\tcp*{Choose one edge $v_iv_j$}
			}
		}
		
		\If{$k<r_k$}{
			\For{$p\leftarrow k+1$ \KwTo $r_k$}{
				$\mathrm{best}\leftarrow \max\{\mathrm{best},\,2+F[r_p+1]\}$\tcp*{Choose two edges $v_iv_k$ and $v_kv_p$}
			}
		}
		
		$F[i]\leftarrow \mathrm{best}$\;
	}
	
	\Return{$F[1]$}\;
\end{algorithm}

Although Algorithm~\ref{EOP_PIG_Algo} is stated only for computing the value $\rho_e^o(G)$, a maximum edge open packing set can also be obtained by a standard backtracking modification. During the computation of $F[i]$, whenever the value of $\mathrm{best}$ is updated, we store the corresponding choice: either the vertex $v_i$ is skipped, one edge $v_iv_j$ is selected, or two edges $v_iv_k$ and $v_kv_p$ are selected. After the table has been filled, we start from $i=1$ and follow these stored choices through the suffix indices $i+1$, $r_j+1$, or $r_p+1$, respectively. The edges recorded along this path form a maximum edge open packing set of $G$. This reconstruction does not change the asymptotic running time of the algorithm.

\subsection{Running time analysis}
The algorithm computes the BCO $\sigma=(v_1,v_2,\ldots,v_n)$ only once. After this ordering is fixed, the values $r_i=r(v_i)$, for $i\in[n]$, can be computed by scanning the adjacency lists according to the ordering. The dynamic programming table has one entry $F[i]$ for each suffix graph $G_i=G[\{v_i,v_{i+1},\ldots,v_n\}]$, together with the terminal value $F[n+1]=0$ for the empty suffix.

The algorithm fills the table from right to left. For each fixed $i$, the first inner loop considers all choices of one edge $v_iv_j$, where $i<j\leq r_i$. This corresponds to the value $\rho'(G_i)$. The second inner loop considers the two-edge case. By the monotonicity of rightmost neighbours in a BCO of a proper interval graph, the maximum over all pairs $i<j\leq r_i$ and $r_i<p\leq r_j$ can be obtained by considering $j=r_i$. Therefore, it is enough to consider $p$ with $r_i<p\leq r_{r_i}$, which is exactly the second loop. This corresponds to the value $\rho''(G_i)$.

Each inner loop has length at most $n$, and there are $n$ choices of $i$. Hence the total time spent in the dynamic programming loops is $O(n^2)$. The table requires $O(n)$ memory. Thus, after the BCO and the values $r_i$ are available, Algorithm~\ref{EOP_PIG_Algo} runs in $O(n^2)$ time. In particular, the algorithm is polynomial. The terminal cases are handled by the initialization $F[n+1]=0$; for example, if $n=1$, the algorithm returns $0$, and if $G\cong K_2$, the two-edge case is skipped and the algorithm returns $1$.

Based on the running time analysis of the algorithm MaxEOP\_PIG(G)\ref{EOP_PIG_Algo}, we have the following theorem.

\begin{theorem}
	Maximum edge open packing can be solved in $O(n^2)$ time for proper interval graphs.
\end{theorem}

\section{EOP set in block graphs}\label{EOP_Sec_3}

In this section, we design a polynomial-time algorithm for computing the edge open packing number of a block graph. Our approach is adapted from the rooted-tree dynamic-programming method for edge open packing in trees given in~\cite{brevsar2024edge}. In a tree, each branch below a vertex is attached through a single edge, and therefore the local decisions around a vertex can be handled directly through its children. For block graphs, the situation is more involved: a cut vertex is connected to its descendants through blocks, and every block is a clique. Hence, an edge chosen inside one block can restrict the choices in several descendant subgraphs at once. Therefore, we cannot simply carry over the tree dynamic program; the states must be refined to capture these block-level interactions. We keep the same general philosophy of recording the status of the distinguished cut vertex, but we introduce block-level transitions that distinguish whether the cut vertex is the centre of a selected star, a leaf of a selected star, not incident with any selected edge, or avoided together with its closed neighbourhood. The main additional work is to combine the contributions of all child cut vertices of a block while respecting the fact that at most one edge from a clique block can belong to an edge open packing set.

Let \(G\) be a connected block graph with blocks \(B_1,\ldots,B_\ell\) and cut vertices \(v_1,v_2,\ldots,v_k\). Each block \(B_i\) is regarded as a subgraph of \(G\), and its vertex set is denoted by \(V(B_i)\). The \emph{cut-tree} of \(G\), denoted by \(T_G\), is defined as follows. For each block \(B_i\), we introduce a vertex \(b_i\) in \(T_G\), called a \emph{block vertex}. The vertex set of \(T_G\) is \(\{b_1,b_2,\ldots,b_\ell\}\cup\{v_1,v_2,\ldots,v_k\}\), and the edge set is \(\{b_i v_j: v_j\in V(B_i), 1\leq i\leq \ell, 1\leq j\leq k\}\). Thus, a block vertex \(b_i\) is adjacent in \(T_G\) exactly to the cut vertices of \(G\) that belong to \(B_i\). The set of all blocks, the set of all cut vertices, and the incidences between blocks and cut vertices can be found in \(O(|V(G)|+|E(G)|)\) time using a depth-first search. Hence, the cut-tree of a block graph can be constructed in \(O(|V(G)|+|E(G)|)\) time~\cite{AhoHopcroftUllman}.

If \(G\) has no cut vertex, then \(G\) is a complete graph. In that case, \(\rho_e^o(G)=0\) when \(|V(G)|=1\), while \(\rho_e^o(G)=1\) when \(|V(G)|\geq 2\). Therefore, in the rest of this section we assume that \(G\) has at least one cut vertex. Fix a cut vertex \(r\) of \(G\), and root the cut-tree \(T_G\) at \(r\). We denote this rooted cut-tree by \(T_G^r\). Throughout this section, all parent-child relations are taken with respect to this fixed rooted tree \(T_G^r\). Thus, the notation \(T_G^r\) is fixed once and for all; the symbol \(G_v\) denotes a subgraph associated with a vertex \(v\), not a new rooted cut-tree.

For an edge set \(S\), let \(V(S)\) denote the set of all endpoints of edges in \(S\). We write \(G[S]\) for the subgraph \(G[V(S)]\). Recall that an edge set \(S\subseteq E(G)\) is an EOP set if every component of \(G[S]\) is a star. In particular, an EOP set of a block graph contains at most one edge from any block, since every block is a clique.

\begin{definition}[Subgraphs determined by the rooted cut-tree]\label{def:block-subgraphs}
	Let \(G\) be a connected block graph with rooted cut-tree \(T_G^r\). For a cut vertex \(v\neq r\), let \(P(v)\) denote the block whose block vertex is the parent of \(v\) in \(T_G^r\). The root \(r\) has no parent block, and in this case we set \(G_r=G\).
	
	For a cut vertex \(v\neq r\), define \(G_v\) to be the connected component containing \(v\) in the graph obtained from \(G\) by deleting all vertices of \(P(v)\setminus\{v\}\). If \(v\) is a non-cut vertex, then \(v\) belongs to a unique block, and we define \(G_v\) to be the one-vertex graph induced by \(v\).
	
	For a cut vertex \(v\), let \(CC(v)\) be the set of blocks \(B\) such that the block vertex corresponding to \(B\) is a child of \(v\) in \(T_G^r\). Similarly, for a block \(B\), let \(CB(B)\) be the set of cut vertices \(u\) such that \(u\) is a child of the block vertex corresponding to \(B\) in \(T_G^r\).
	
	Now let \(v\) be a cut vertex and let \(B\in CC(v)\). We define \(G_v^B\) as the subgraph of \(G_v\) induced by the vertices lying below the child block \(B\), after removing \(v\). Equivalently, \(G_v^B=G_v[(V(B)\setminus\{v\})\cup \bigcup_{u\in CB(B)}V(G_u)]\). Thus, \(G_v^B\) consists of the vertices of the block \(B\) except \(v\), together with all subgraphs rooted at the child cut vertices of \(B\).
\end{definition}

The graphs \(G_v\) and \(G_v^B\) are illustrated in Figure~\ref{fig:block_cut_tree}.

\par\medskip
\noindent
\begin{minipage}{\linewidth}
	\centering
	\includegraphics[width=0.85\textwidth]{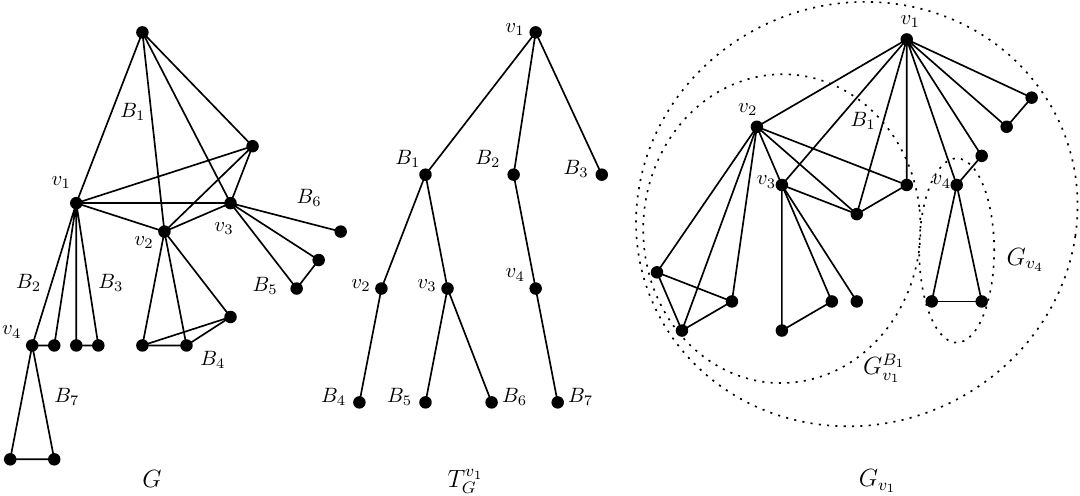}
	\captionof{figure}{A block graph \(G\), its rooted cut-tree \(T_G^{v_1}\), and the associated subgraphs. The graph on the left is the block graph \(G\) in its original drawing, while the drawing on the right represents the same graph arranged with \(v_1\) as the root, so that the subgraphs \(G_{v_1}\) and \(G_{v_1}^{B_1}\) can be displayed more clearly. The rooted cut-tree \(T_G^{v_1}\) is shown in the middle.}
	\label{fig:block_cut_tree}
\end{minipage}
\par\medskip

To design the dynamic program, we use five parameters for each graph \(G_v\) defined above. For a one-vertex graph \(G_v\), all five parameters below are set equal to \(0\). For the remaining cases, we define the following values:
\begin{itemize}
	\item \(\rho_e^o(G_v)\) is the maximum size of an EOP set in \(G_v\).
	\item \(\rho_e^c(G_v)\) is the maximum size of an EOP set \(S\) in \(G_v\) such that \(v\) is the centre of a star in \(G_v[S]\).
	\item \(\rho_e^l(G_v)\) is the maximum size of an EOP set \(S\) in \(G_v\) such that \(v\) is a leaf of a star in \(G_v[S]\).
	\item \(\rho_e^{\prime}(G_v)\) is the maximum size of an EOP set \(S\) in \(G_v\) such that no edge of \(S\) is incident with \(v\).
	\item \(\rho_e^{\prime\prime}(G_v)\) is the maximum size of an EOP set \(S\) in \(G_v\) such that \(N_{G_v}[v]\cap V(S)=\emptyset\).
\end{itemize}
The condition in the definition of \(\rho_e^{\prime}(G_v)\) may also be written as \(v\notin V(S)\), or equivalently as \(v\notin V(G_v[S])\). We avoid writing \(v\notin S\), since \(S\) is a set of edges, not a set of vertices.

For every vertex \(v\) for which \(G_v\) is defined, an EOP set \(S\) of \(G_v\) has exactly one of the following three statuses at \(v\): either \(v\) is the centre of a star in \(G_v[S]\), or \(v\) is a leaf of a star in \(G_v[S]\), or no edge of \(S\) is incident with \(v\). If \(v\) belongs to a \(K_{1,1}\)-component of \(G_v[S]\), then \(v\) may be regarded either as the centre or as a leaf of that component. Hence we obtain the following observation.

\begin{observation}\label{Block_obs_1}
	For every vertex \(v\) for which \(G_v\) is defined, every EOP set of \(G_v\) falls into one of three cases: either \(v\) is the centre of a selected star, \(v\) is a leaf of a selected star, or no selected edge is incident with \(v\). Therefore, \(\rho_e^o(G_v)\) is the maximum of the three values \(\rho_e^c(G_v)\), \(\rho_e^l(G_v)\), and \(\rho_e^{\prime}(G_v)\). In particular, since \(G_r=G\), we have \(\rho_e^o(G)=\rho_e^o(G_r)\), and this is the maximum of \(\rho_e^c(G_r)\), \(\rho_e^l(G_r)\), and \(\rho_e^{\prime}(G_r)\).
\end{observation}

We now give the recurrence relations. All sums over empty sets are interpreted as \(0\). The values for a non-cut vertex are also taken as \(0\), because the corresponding graph \(G_v\) has only one vertex.

For a cut vertex \(v\) and a child block \(B\in CC(v)\), set \(R_B=\sum_{w\in CB(B)}\rho_e^{\prime}(G_w)\). This is the value obtained from the descendants of \(B\) when no vertex of \(B\) is used as an endpoint of a selected edge.

First we need the value of \(\rho_e^o(G_v^B)\), because this value appears when no selected edge is incident with \(v\). Let \(C_B=V(B)\setminus\{v\}\). Define the following three quantities:
\begin{equation*}
	\Theta_0(v,B)=R_B .
\end{equation*}
\begin{equation*}
	\begin{aligned}
		\Theta_1(v,B)
		=\max_{x\in C_B}\Bigl\{&
		\max\{\rho_e^c(G_x),\rho_e^l(G_x)\}+R_B\\
		&-\rho_e^{\prime}(G_x)\Bigr\}.
	\end{aligned}
\end{equation*}
\begin{equation*}
	\begin{aligned}
		\Theta_2(v,B)
		=\max_{\substack{x,y\in C_B\\ x\neq y}}\Bigl\{&
		1+\max\{\rho_e^c(G_x),\rho_e^{\prime\prime}(G_x)\}\\
		&+\rho_e^{\prime\prime}(G_y)+R_B\\
		&-\rho_e^{\prime}(G_x)-\rho_e^{\prime}(G_y)\Bigr\}.
	\end{aligned}
\end{equation*}
If the index set of \(\Theta_1(v,B)\) or \(\Theta_2(v,B)\) is empty, then the corresponding value is taken as \(0\). Finally, let \(\Theta(v,B)=\max\{\Theta_0(v,B),\Theta_1(v,B),\Theta_2(v,B)\}\).

\begin{lemma}\label{Block_lemma_branch}
	Let \(G\) be a connected block graph with rooted cut-tree \(T_G^r\). If \(v\) is a cut vertex and \(B\in CC(v)\), then \(\rho_e^o(G_v^B)=\Theta(v,B)\).
\end{lemma}

\begin{proof}
	Let \(S\) be an EOP set of \(G_v^B\), and consider the set of vertices of the clique \(C_B=V(B)\setminus\{v\}\) that appear as endpoints of edges in \(S\). Since \(C_B\) is a clique and every component of \(G_v^B[S]\) must be a star, at most two vertices of \(C_B\) can be endpoints.
	
	If no vertex of \(C_B\) is an endpoint, then each child branch \(G_w\), where \(w\in CB(B)\), contributes at most \(\rho_e^{\prime}(G_w)\). This gives the upper bound \(\Theta_0(v,B)\).
	
	Suppose exactly one vertex \(x\in C_B\) is an endpoint. Then no edge of the top clique \(C_B\) is selected. The endpoint \(x\) must be used inside its own branch, and hence the contribution of \(G_x\) is at most \(\max\{\rho_e^c(G_x),\rho_e^l(G_x)\}\). Every other child branch contributes under the condition that its root is not incident with a selected edge. This gives the upper bound \(\Theta_1(v,B)\).
	
	Finally, suppose two vertices \(x,y\in C_B\) are endpoints. Since \(C_B\) is a clique, the edge \(xy\) is present. In an optimal EOP set with these two top endpoints, the edge \(xy\) must be the selected edge from the top block; otherwise the edge \(xy\) would join two selected components and the induced subgraph would not be a star. One of \(x\) and \(y\) acts as the centre of the resulting star and the other as a leaf. If \(x\) is the centre and \(y\) is the leaf, then the branch below \(x\) contributes at most \(\max\{\rho_e^c(G_x),\rho_e^{\prime\prime}(G_x)\}\), while the branch below \(y\) contributes at most \(\rho_e^{\prime\prime}(G_y)\). All other child branches contribute through their \(\rho_e^{\prime}\)-values. This gives the upper bound \(\Theta_2(v,B)\). The opposite orientation is already covered by interchanging \(x\) and \(y\) in the maximum.
	
	The reverse inequality follows by choosing optimal sets that realize the corresponding terms in \(\Theta_0(v,B)\), \(\Theta_1(v,B)\), and \(\Theta_2(v,B)\). In each case, the selected edges inside different child branches are compatible because the only possible interaction between them is through the clique \(C_B\), and the three cases above control exactly how many vertices of \(C_B\) are used as endpoints. Therefore, the maximum of these three quantities is attainable, and the lemma follows.
\end{proof}

We next compute the four restricted values for \(G_v\), where \(v\) is a cut vertex. For a block \(B\in CC(v)\), define
\begin{equation*}
	\begin{aligned}
		\mathcal A_B
		=\max_{u\in V(B)\setminus\{v\}}\Bigl\{&
		1+\rho_e^{\prime\prime}(G_u)\\
		&+\sum_{w\in CB(B)\setminus\{u\}}\rho_e^{\prime}(G_w)\Bigr\}.
	\end{aligned}
\end{equation*}
\begin{equation*}
	\overline{\mathcal A}_B=\sum_{w\in CB(B)}\rho_e^{\prime}(G_w).
\end{equation*}
The value \(\mathcal A_B\) corresponds to the case where the selected set contains one edge \(vu\) from the block \(B\), while \(\overline{\mathcal A}_B\) corresponds to the case where no selected edge from \(B\) is incident with \(v\). We call \(B\) a \emph{Type 1} child block of \(v\) if \(\mathcal A_B\geq \overline{\mathcal A}_B\), and a \emph{Type 2} child block otherwise.

\begin{lemma}\label{Block_lemma_1}
	Let \(G\) be a connected block graph with rooted cut-tree \(T_G^r\), and let \(v\) be a cut vertex. If \(CC(v)=\emptyset\), then \(\rho_e^c(G_v)=0\). Otherwise, \(\rho_e^c(G_v)\) is computed as follows. If \(CC(v)\) has a Type 1 block, then
	\begin{equation*}
		\begin{aligned}
			\rho_e^c(G_v)=&
			\sum_{\substack{B\in CC(v)\\ B\text{ Type 1}}}\mathcal A_B\\
			&+\sum_{\substack{B\in CC(v)\\ B\text{ Type 2}}}\overline{\mathcal A}_B .
		\end{aligned}
	\end{equation*}
	If every block in \(CC(v)\) is Type 2, then
	\begin{equation*}
		\rho_e^c(G_v)=
		\sum_{B\in CC(v)}\overline{\mathcal A}_B
		+\max_{B\in CC(v)}\{\mathcal A_B-\overline{\mathcal A}_B\}.
	\end{equation*}
\end{lemma}

\begin{proof}
	Let \(S\) be an EOP set of \(G_v\) such that \(v\) is the centre of a star in \(G_v[S]\). Then at least one child block of \(v\) contains a selected edge incident with \(v\). Since each child block is a clique, at most one edge can be selected from such a block. If the selected edge in \(B\) is \(vu\), then no vertex in \(N_{G_u}[u]\) can be an endpoint of an edge selected inside \(G_u\), and hence the branch below \(u\) contributes at most \(\rho_e^{\prime\prime}(G_u)\). Every other child cut vertex \(w\in CB(B)\setminus\{u\}\) must not be incident with a selected edge, and hence its branch contributes at most \(\rho_e^{\prime}(G_w)\). Therefore, a block that contributes an edge incident with \(v\) contributes at most \(\mathcal A_B\). A block that contributes no such edge contributes at most \(\overline{\mathcal A}_B\).
	
	Thus, if there is a Type 1 block, the best choice is to take an edge incident with \(v\) from every Type 1 block and no such edge from every Type 2 block. These selected edges have the common endpoint \(v\), and their other endpoints lie in distinct child blocks, so they form a star centered at \(v\). This gives the first expression. If all child blocks are Type 2, then we must still choose one child block to provide an edge incident with \(v\), and the best choice is the block maximizing \(\mathcal A_B-\overline{\mathcal A}_B\). This gives the second expression. The constructions described above also attain the corresponding values, so the equality follows.
\end{proof}

For the case where \(v\) is a leaf of a selected star, exactly one child block of \(v\) supplies the edge incident with \(v\). For \(B\in CC(v)\), define
\begin{equation*}
	\begin{aligned}
		\mathcal M_B
		=\max_{u\in V(B)\setminus\{v\}}\Bigl\{&
		\max\{\rho_e^c(G_u),\rho_e^{\prime\prime}(G_u)\}+1\\
		&+\sum_{w\in CB(B)\setminus\{u\}}\rho_e^{\prime}(G_w)\Bigr\}.
	\end{aligned}
\end{equation*}
\begin{equation*}
	\overline{\mathcal M}_B
	=\sum_{B'\in CC(v)\setminus\{B\}}\sum_{w\in CB(B')}\rho_e^{\prime}(G_w).
\end{equation*}

\begin{lemma}\label{Block_lemma_2}
	Let \(G\) be a connected block graph with rooted cut-tree \(T_G^r\). If \(v\) is a cut vertex, then \(\rho_e^l(G_v)=\displaystyle\max_{B\in CC(v)}\{\mathcal M_B+\overline{\mathcal M}_B\}\). If \(CC(v)=\emptyset\), then \(\rho_e^l(G_v)=0\).
\end{lemma}

\begin{proof}
	Let \(S\) be an EOP set of \(G_v\) such that \(v\) is a leaf of a star in \(G_v[S]\). Then there is exactly one child block \(B\in CC(v)\) containing the selected edge incident with \(v\). Let this edge be \(vu\), where \(u\in V(B)\setminus\{v\}\). In the resulting star, the vertex \(u\) is the centre. Inside the branch \(G_u\), either \(u\) is already the centre of a selected star, or no vertex in \(N_{G_u}[u]\) is used by the selected set. These two possibilities contribute at most \(\rho_e^c(G_u)\) and \(\rho_e^{\prime\prime}(G_u)\), respectively. Every other child cut vertex \(w\in CB(B)\setminus\{u\}\) cannot be incident with a selected edge, so its branch contributes at most \(\rho_e^{\prime}(G_w)\). This gives \(\mathcal M_B\).
	
	For every other child block \(B'\neq B\), no selected edge can be incident with \(v\). Hence each child cut vertex \(w\in CB(B')\) contributes through \(\rho_e^{\prime}(G_w)\), giving \(\overline{\mathcal M}_B\). Maximizing over the unique block \(B\) that contains the edge incident with \(v\) gives the stated formula. The reverse inequality follows by taking optimal sets for the chosen block and for all remaining branches; these sets are compatible because the only selected edge incident with \(v\) is \(vu\), and the construction makes \(v\) a leaf of the resulting star.
\end{proof}

\begin{lemma}\label{Block_lemma_3}
	Let \(G\) be a connected block graph with rooted cut-tree \(T_G^r\), and let \(v\) be a cut vertex of \(G\). Then \(\rho_e^{\prime}(G_v)\) is obtained by summing the optimum values over all branches below the child blocks of \(v\). More precisely, \(\rho_e^{\prime}(G_v)=\sum_{B\in CC(v)}\rho_e^o(G_v^B)\). By Lemma~\ref{Block_lemma_branch}, this is equivalently \(\rho_e^{\prime}(G_v)=\sum_{B\in CC(v)}\Theta(v,B)\).
\end{lemma}

\begin{proof}
	Let \(S\) be an EOP set of \(G_v\) such that no edge of \(S\) is incident with \(v\). Then the selected edges lie completely inside the subgraphs \(G_v^B\), where \(B\in CC(v)\). For each such block, let \(S_B=S\cap E(G_v^B)\). Then \(S_B\) is an EOP set of \(G_v^B\), and hence \(|S_B|\leq \rho_e^o(G_v^B)\). Since the graphs \(G_v^B\) are edge-disjoint and interact only through the vertex \(v\), which is not incident with any selected edge, we obtain \(|S|\leq \sum_{B\in CC(v)}\rho_e^o(G_v^B)\).
	
	Conversely, for each \(B\in CC(v)\), choose a maximum EOP set of \(G_v^B\). The union of these sets is an EOP set of \(G_v\) in which no edge is incident with \(v\). Therefore, \(\rho_e^{\prime}(G_v)=\sum_{B\in CC(v)}\rho_e^o(G_v^B)\). By Lemma~\ref{Block_lemma_branch}, this is equal to \(\sum_{B\in CC(v)}\Theta(v,B)\).
\end{proof}

\begin{lemma}\label{Block_lemma_4}
	Let \(G\) be a connected block graph with rooted cut-tree \(T_G^r\), and let \(v\) be a cut vertex of \(G\). Then \(\rho_e^{\prime\prime}(G_v)\) is obtained by summing the values \(\rho_e^{\prime}(G_u)\) over all child cut vertices \(u\) lying below the child blocks of \(v\). More precisely, \(\rho_e^{\prime\prime}(G_v)=\sum_{B\in CC(v)}\sum_{u\in CB(B)}\rho_e^{\prime}(G_u)\).
\end{lemma}

\begin{proof}
	Let \(S\) be an EOP set of \(G_v\) such that \(N_{G_v}[v]\cap V(S)=\emptyset\). Since every vertex of \(V(B)\setminus\{v\}\) is adjacent to \(v\) for each \(B\in CC(v)\), no edge of \(S\) can have an endpoint in \(V(B)\setminus\{v\}\). Therefore, the only possible selected edges lie inside the descendant subgraphs \(G_u\), where \(u\in CB(B)\) and \(B\in CC(v)\). Moreover, no such selected edge can be incident with \(u\), because \(u\) is adjacent to \(v\). Thus, the contribution from each \(G_u\) is at most \(\rho_e^{\prime}(G_u)\), and we get \(|S|\leq \sum_{B\in CC(v)}\sum_{u\in CB(B)}\rho_e^{\prime}(G_u)\).
	
	Conversely, for every \(B\in CC(v)\) and every \(u\in CB(B)\), choose an EOP set in \(G_u\) realizing \(\rho_e^{\prime}(G_u)\). The union of all these sets is an EOP set of \(G_v\), and no endpoint of a selected edge lies in \(N_{G_v}[v]\). Hence the upper bound is attained.
\end{proof}

The preceding observation and lemmas lead to the following dynamic program. The cut vertices are processed in non-increasing order of their distance from the root \(r\) in \(T_G^r\), so every child cut vertex is processed before its parent.

\begin{algorithm}[h]
	\caption{\textsc{MaxEOP-Block}\((G)\)}\label{EOP_BG_Algo}
	\KwIn{A connected block graph \(G=(V,E)\).}
	\KwOut{\(\rho_e^o(G)\), the edge open packing number of \(G\).}
	\If{\(G\) has no cut vertex}{
		\If{\(|V(G)|=1\)}{\Return \(0\)\;}
		\Else{\Return \(1\)\;}
	}
	Choose a cut vertex \(r\), construct the cut-tree \(T_G\), and root it at \(r\), obtaining \(T_G^r\)\;
	For every non-cut vertex \(x\), set \(\rho_e^o(G_x)=\rho_e^c(G_x)=\rho_e^l(G_x)=\rho_e^{\prime}(G_x)=\rho_e^{\prime\prime}(G_x)=0\)\;
	Let \(\sigma=(v_1,v_2,\ldots,v_k)\) be the cut vertices ordered so that \(d_{T_G^r}(r,v_1)\geq d_{T_G^r}(r,v_2)\geq \cdots\geq d_{T_G^r}(r,v_k)\)\;
	\For{\(i=1\) to \(k\)}{
		\For{each \(B\in CC(v_i)\)}{
			Compute \(R_B=\sum_{w\in CB(B)}\rho_e^{\prime}(G_w)\)\;
			Compute \(\Theta(v_i,B)\) using Lemma~\ref{Block_lemma_branch}\;
			Compute \(\mathcal A_B\), \(\overline{\mathcal A}_B\), \(\mathcal M_B\), and \(\overline{\mathcal M}_B\)\;
			Mark \(B\) as Type 1 if \(\mathcal A_B\geq \overline{\mathcal A}_B\), and as Type 2 otherwise\;
		}
		Set \(\rho_e^{\prime\prime}(G_{v_i})=\sum_{B\in CC(v_i)}\sum_{u\in CB(B)}\rho_e^{\prime}(G_u)\)\;
		Set \(\rho_e^{\prime}(G_{v_i})=\sum_{B\in CC(v_i)}\Theta(v_i,B)\)\;
		\uIf{\(CC(v_i)=\emptyset\)}{
			Set \(\rho_e^c(G_{v_i})=0\) and \(\rho_e^l(G_{v_i})=0\)\;
		}
		\uElseIf{every block in \(CC(v_i)\) is Type 2}{
			Set \(\rho_e^c(G_{v_i})=\sum_{B\in CC(v_i)}\overline{\mathcal A}_B+\max_{B\in CC(v_i)}\{\mathcal A_B-\overline{\mathcal A}_B\}\)\;
			Set \(\rho_e^l(G_{v_i})=\max_{B\in CC(v_i)}\{\mathcal M_B+\overline{\mathcal M}_B\}\)\;
		}
		\Else{
			Set \(\rho_e^c(G_{v_i})=\sum_{\substack{B\in CC(v_i)\\ B\text{ Type 1}}}\mathcal A_B+\sum_{\substack{B\in CC(v_i)\\ B\text{ Type 2}}}\overline{\mathcal A}_B\)\;
			Set \(\rho_e^l(G_{v_i})=\max_{B\in CC(v_i)}\{\mathcal M_B+\overline{\mathcal M}_B\}\)\;
		}
		Set \(\rho_e^o(G_{v_i})=\max\{\rho_e^c(G_{v_i}),\rho_e^l(G_{v_i}),\rho_e^{\prime}(G_{v_i})\}\)\;
	}
	\Return \(\rho_e^o(G_r)\)\;
\end{algorithm}

\begin{theorem}\label{thm:block-eop-correctness}
	Algorithm~\ref{EOP_BG_Algo} correctly computes \(\rho_e^o(G)\) for every connected block graph \(G\).
\end{theorem}

\begin{proof}
	If \(G\) has no cut vertex, then \(G\) is complete, and the value returned by the algorithm is immediate. Otherwise, the algorithm fixes a rooted cut-tree \(T_G^r\) and processes the cut vertices from the leaves towards the root. Hence, when a cut vertex \(v_i\) is processed, all values associated with the child cut vertices below \(v_i\) have already been computed.
	
	For each child block \(B\in CC(v_i)\), Lemma~\ref{Block_lemma_branch} gives the correct value of \(\rho_e^o(G_{v_i}^B)\). Lemmas~\ref{Block_lemma_1}, \ref{Block_lemma_2}, \ref{Block_lemma_3}, and \ref{Block_lemma_4} then give the correct values of \(\rho_e^c(G_{v_i})\), \(\rho_e^l(G_{v_i})\), \(\rho_e^{\prime}(G_{v_i})\), and \(\rho_e^{\prime\prime}(G_{v_i})\), respectively. Finally, Observation~\ref{Block_obs_1} gives \(\rho_e^o(G_{v_i})\). By induction over the processing order, every value computed by the algorithm is correct. Since \(G_r=G\), the returned value \(\rho_e^o(G_r)\) is exactly \(\rho_e^o(G)\).
\end{proof}

\subsection{Running time}\label{subsec:block-running-time}

Let \(n=|V(G)|\) and \(m=|E(G)|\). The blocks, cut vertices, and the cut-tree \(T_G\) can be computed in \(O(n+m)\) time. Rooting \(T_G\) and producing the postorder of the cut vertices also takes linear time in the size of \(T_G\), which is \(O(n)\).

For a fixed child block \(B\), the quantities \(R_B\), \(\mathcal A_B\), \(\overline{\mathcal A}_B\), \(\mathcal M_B\), and \(\overline{\mathcal M}_B\) can be computed by scanning the vertices of \(B\) and the child cut vertices in \(CB(B)\). The only term that may require pairs of vertices is \(\Theta_2(v,B)\), and it can be computed in \(O(|V(B)|^2)\) time by checking all ordered pairs \((x,y)\) with \(x,y\in V(B)\setminus\{v\}\) and \(x\neq y\). Therefore, the total time spent over all blocks is \(O(\sum_B |V(B)|^2)\). Since the blocks of a block graph are edge-disjoint cliques, \(\sum_B \binom{|V(B)|}{2}=m\), and hence \(\sum_B |V(B)|^2=O(n+m)\). Thus the total running time of Algorithm~\ref{EOP_BG_Algo} is \(O(n+m)\). In particular, since a block graph may have \(m=O(n^2)\), the algorithm runs in \(O(n^2)\) time in terms of the number of vertices alone. The space usage is \(O(n)\), apart from the storage of the input graph.

\begin{theorem}\label{thm:block-eop-polynomial}
	Maximum edge open packing can be solved in polynomial time for block graphs. More precisely, for a block graph given by adjacency lists, Algorithm~\ref{EOP_BG_Algo} runs in \(O(n+m)\) time.
\end{theorem}

If the input block graph is disconnected, then the algorithm is applied independently to each connected component, and the final value is the sum of the values over all components.


\section{EOP set in split graphs}\label{EOP_Sec_4}

In this section, we show that the edge open packing number of a split graph can be found in polynomial time. Throughout this section, let \(G=(S\cup K,E)\) be a split graph, where \(S\) is the independent set and \(K\) is the clique of \(G\). We assume that \(K\) is a maximum clique of \(G\), that is, \(\omega(G)=|K|\). For \(x\in K\), let \(N_S(x)=N(x)\cap S\) and let \(l_x=|N_S(x)|\). For distinct vertices \(x,y\in K\), let \(l_{xy}=|N_S(x)\cap N_S(y)|\). If \(|K|\geq 2\), define \(\rho_1(G)=1+\max\{l_x-l_{xy}:x,y\in K,\ x\neq y\}\), and if \(|K|\leq 1\), define \(\rho_1(G)=0\). Also, define \(\rho_2(G)=\max\{l_x:x\in K\}\), where this maximum is taken to be \(0\) when \(K=\emptyset\).

\begin{theorem}\label{thm:eop-split}
	Let \(G=(S\cup K,E)\) be a split graph. Then \(\rho_e^o(G)=\max\{\rho_1(G),\rho_2(G)\}\).
\end{theorem}

\begin{proof}
	We prove the two inequalities separately. First, we show that \(\rho_e^o(G)\geq \max\{\rho_1(G),\rho_2(G)\}\). For any \(x\in K\), the set \(D_x=\{xs:s\in N_S(x)\}\) is an EOP set, since all its edges form a star centred at \(x\) and \(S\) is independent. Hence, \(\rho_e^o(G)\geq l_x\) for every \(x\in K\), and therefore, \(\rho_e^o(G)\geq \rho_2(G)\).
	
	Now let \(x,y\in K\) be distinct. Consider the set \(D_{x,y}=\{xy\}\cup \{xs:s\in N_S(x)\setminus N_S(y)\}\). This is also an EOP set. Indeed, all edges of \(D_{x,y}\) form a star centred at \(x\), and no vertex \(s\in N_S(x)\setminus N_S(y)\) is adjacent to \(y\). Thus no edge \(ys\) can create a common edge between \(xy\) and \(xs\). Hence, \(|D_{x,y}|=1+l_x-l_{xy}\). Taking the maximum over all ordered pairs \(x,y\in K\) with \(x\neq y\), we get \(\rho_e^o(G)\geq \rho_1(G)\). Therefore, \(\rho_e^o(G)\geq \max\{\rho_1(G),\rho_2(G)\}\).
	
	Conversely, let \(D\) be a maximum EOP set of \(G\). We prove that \(|D|\leq \max\{\rho_1(G),\rho_2(G)\}\). Since \(K\) is a clique, \(D\) contains at most one edge from \(E(G[K])\). First suppose that \(D\cap E(G[K])=\emptyset\). Then no two edges of \(D\) can be incident with two different vertices of \(K\), because the edge between those two clique vertices would be a common edge. Hence, all edges of \(D\) are incident with a single vertex \(x\in K\), and their other endvertices lie in \(S\). Therefore, \(|D|\leq l_x\leq \rho_2(G)\).
	
	Now suppose that \(D\) contains an edge \(xy\in E(G[K])\). Since \(K\) is a clique, no edge of \(D\setminus\{xy\}\) can be incident with a vertex of \(K\setminus\{x,y\}\). Moreover, the remaining edges of \(D\) cannot be incident with both \(x\) and \(y\), because then the edge \(xy\) would be a common edge between two selected edges. Thus all edges of \(D\setminus\{xy\}\), if any, are incident with exactly one of \(x\) and \(y\). If they are incident with \(x\), then their endvertices in \(S\) cannot be adjacent to \(y\); otherwise, for a vertex \(s\in N_S(x)\cap N_S(y)\), the edge \(ys\) would be a common edge between \(xy\) and \(xs\). Hence, in this case \(|D|\leq 1+l_x-l_{xy}\leq \rho_1(G)\). Similarly, if the remaining edges are incident with \(y\), then \(|D|\leq 1+l_y-l_{xy}\leq \rho_1(G)\). Therefore, in every case with \(xy\in D\), we have \(|D|\leq \rho_1(G)\).
	
	Combining the two cases, every maximum EOP set has cardinality at most \(\max\{\rho_1(G),\rho_2(G)\}\). Together with the lower bound proved above, this gives us \(\rho_e^o(G)=\max\{\rho_1(G),\rho_2(G)\}\).
\end{proof}
%
\subsection{Running time analysis}

Let \(n=|V(G)|\). If the split partition \(S\cup K\) is not given as part of the input, then we first find such a partition, with \(K\) chosen as a maximum clique, using a standard linear-time recognition algorithm for split graphs. Since \(|E(G)|\leq n^2\), this step takes \(O(n^2)\) time in terms of \(n\).

Once the partition is available, the values \(l_x=|N(x)\cap S|\), for all \(x\in K\), can be computed by scanning the edges between \(S\) and \(K\). This takes \(O(n^2)\) time. To compute \(l_{xy}=|N(x)\cap N(y)\cap S|\) for all distinct \(x,y\in K\), we initialize a table indexed by unordered pairs of vertices of \(K\). Then, for each vertex \(s\in S\), we scan \(N(s)\cap K\) and increase the table entry for every pair of vertices in \(N(s)\cap K\). The time needed for this step is \(O(\sum_{s\in S}d_K(s)^2)\), where \(d_K(s)=|N(s)\cap K|\). Since \(d_K(s)\leq n\) and \(|S|\leq n\), we have \(\sum_{s\in S}d_K(s)^2\leq n^3\). After these values are computed, \(\rho_1(G)\) is obtained by scanning all pairs of vertices of \(K\), and \(\rho_2(G)\) is obtained by scanning all vertices of \(K\), which takes at most \(O(n^2)\) additional time. Therefore, the total running time is \(O(n^3)\).

\begin{theorem}\label{thm:eop-split-polytime}
	The edge open packing number of a split graph on \(n\) vertices can be found in \(O(n^3)\) time.
\end{theorem}

\section{Conclusion}\label{EOP_Sec_5}

In this paper, we have partially addressed the open question posed by Brešar and Samadi concerning the complexity of computing the edge open packing number in chordal graphs. We studied three important subclasses of chordal graphs. First, we designed an \(O(n^2)\)-time algorithm for proper interval graphs. Next, we presented an \(O(n+m)\)-time algorithm for block graphs, where \(n=|V(G)|\) and \(m=|E(G)|\). Finally, we showed that the edge open packing number of a split graph on \(n\) vertices can be found in \(O(n^3)\) time.

These results give further evidence that the edge open packing problem may admit efficient algorithms on several structured subclasses of chordal graphs. It would be interesting to investigate whether the problem can be solved in polynomial time for all chordal graphs. Other open problems mentioned in~\cite{chelladurai2022edge} and~\cite{brevsar2024edge} also remain natural directions for future research.

\bibliographystyle{plain}
\bibliography{EOP_chordal_Bib}

\end{document}